\documentclass[11pt,reqno]{amsart}

\usepackage{amsmath, amsfonts, amsthm, amssymb, color, cite}

\newtheorem{Theorem}{Theorem}[section]
\newtheorem{Lemma}{Lemma}[section]

\theoremstyle{definition}
\newtheorem{Definition}{Definition}[section]

\theoremstyle{remark}
\newtheorem{Remark}{Remark}[section]

\numberwithin{equation}{section}
\allowdisplaybreaks

\newcommand{\R}{{\mathbb R}}

\newcommand{\eps}{ \varepsilon}
\def\f{\frac}

\def\hf1{^\f{1}{1-\xi^2}}

\def\be{\begin{equation}}
\def\en{\end{equation}}
\def\bs{\begin{split}}
\def\es{\end{split}}
\def\ba{\begin{align}}
\def\ea{\end{align}}

\author[D.F.Yuan]{Difan Yuan}
\address{School of Mathematical Sciences,  Beijing Normal University and Laboratory of Mathematics and Complex Systems, Ministry of Education, Beijing 100875, China,}
\email{yuandf@amss.ac.cn}

\title[Cylindrical symmetry]
{Global solutions to rotating motion of isentropic flows with cylindrical symmetry}


\date{\today}

\begin{document}
\begin{abstract}
We are concerned with global weak solutions to the isentropic compressible Euler equations with cylindrically symmetric rotating structure, in which the origin is included. Due to the presence of the singularity at the origin, only the case excluding the origin $|\vec{x}|\geq1$ has been considered by Chen-Glimm \cite{Chen3}. The convergence and consistency of the approximate solutions are proved by using $L^{\infty}$ compensated compactness framework and vanishing viscosity method. We observe that if the blast wave initially moves outwards and if initial density and velocity decay to zero at certain algebraic rate near the origin, then the density and velocity decay at the same rate for any positive time. In particular, the initial normal velocity is assumed to be non-negative, and there is no restriction on the sign of initial angular velocity.
\end{abstract}

\maketitle
\section{Introduction}
%
In this paper, we prove the existence result for the global entropy solutions to two dimensional isentropic compressible Euler equations with cylindrically symmetric rotating flow including the origin. The compressible Euler equations are of the following conservative form:
\begin{eqnarray}\label{cylindricalsym}
\left\{ \begin{array}{ll}
\displaystyle \rho_t+\nabla\cdot \vec{m}=0,\\
\displaystyle \vec{m}_t+\nabla\cdot\left(\frac{\vec{m}\otimes\vec{m}}{\rho}\right)+\nabla p=0.
\end{array}
\right.
\end{eqnarray}
We are interested in cylindrically symmetric flow to system \eqref{cylindricalsym} with the form:\\
\begin{equation}\label{transform}
(\rho,\vec{m})(\vec{x},t)=\left(\rho(x,t),m(x,t)\frac{\vec{x}}{x}+\tilde{m}(x,t)\frac{(-x_2,x_1)}{x}\right), x=|\vec{x}|.
\end{equation}
Then $(\rho(x,t),m(x,t),\tilde{m}(x,t))$ in \eqref{transform} is governed by the one-dimensional Euler equations with geometric source terms:
\begin{eqnarray}\label{geometric}
\left\{ \begin{array}{ll}
\displaystyle \rho_t+m_x=-\frac{1}{x}m,x\in[0,+\infty),t>0,\\
\displaystyle m_t+\left(\frac{m^2}{\rho}+p(\rho)\right)_x=-\frac{1}{x}\frac{m^{2}-\tilde{m}^2}{\rho},x\in[0,+\infty),t>0,\\
\displaystyle
\tilde{m}_t+\left(\frac{m\tilde{m}}{\rho}\right)_x=-\frac{2}{x}\frac{m\tilde{m}}{\rho},
x\in[0,+\infty),t>0.\\
\end{array}
\right.
\end{eqnarray}
Here $\rho, m,\tilde{ m}$ and $p(\rho)$ represent the density, the normal and angular momentum and the pressure of the gas separately. For polytropic gas, $p(\rho)=p_0\rho^\gamma,$ with $p_0=\frac{\theta^{2}}{\gamma},\theta=\frac{\gamma-1}{2}$ and adiabatic exponent $\gamma>1.$

Various studies on spherical or cylindrical symmetric flow has been aroused attention in recent decades, motivated by many important physical phenomena, such as the stellar dynamics including supernova formation, inertial confinement fusion.  See \cite{Courant,Dafermos,Drake} and references therein.

Consider the initial value problem for \eqref{geometric} with geometric structure:
\begin{eqnarray}\label{geometricc2}
\left\{ \begin{array}{ll}
\displaystyle v_t+F(v)_x=G(x,v),x\in[0,+\infty),t\in[0,+\infty),\\
\displaystyle v|_{t=0}=v_0(x), x\in[0,+\infty),\\
\end{array}
\right.
\end{eqnarray}
where the vector
$v=(\rho, m,\tilde{m})^\top,$ the flux $F(v)=(m, \frac{m^2}{\rho}+p(\rho),\frac{m\tilde{m}}{\rho})^\top,$ the source term
$G(x,v)=(a(x)m,a(x)\frac{m^{2}-\tilde{m}^{2}}{\rho},
a(x)\frac{2m\tilde{m}}{\rho})^\top, $with $a(x)=-\frac{A'(x)}{A(x)}=-\frac{1}{x}, A(x)=2\pi x.$

Here, we consider the Cauchy problem of the compressible Euler equation \eqref{cylindricalsym} with cylindrical symmetric initial data
\begin{equation}\label{sphericalinitial2}
\begin{aligned}
&(\rho, m,\tilde{m})|_{t=0}=(\rho_0(x), m_0(x),\tilde{m}_0(x)),x\geq0.\\
\end{aligned}
\end{equation}

There have been extensive studies on one dimensional isentropic gas dynamics. The first global existence of weak entropy solution with large initial data was established in Diperna \cite{Diperna} by introducing vanishing viscosity method for $\gamma=1+\frac{2}{2n+1},$ where $n\geq2 $ is any positive integer. Ding-Chen-Luo \cite{Ding} and Chen \cite{chen} extended the result for general values $\gamma\in(1, \frac{5}{3}]$ by using Lax-Friedrichs scheme. Lions-Perthame-Tadmor \cite{Lions1} and Lions-Perthame-Souganidis \cite{Lions2} dealt with the case $\gamma>\frac{5}{3}.$ When $\gamma=1,$ Huang-Wang \cite{HuangWang} obtained global existence of $L^{\infty}$ entropy solutions of isothermal gases by adopting analytic extension method.

 For the spherical symmetry solutions of compressible Euler equations, Makino et al. studied the global weak entropy solution outside a solid ball for $\gamma=1.$ The global existence for general $\gamma$ was first studied by Chen-Glimm \cite{Chen2}, and then by Tsuge  \cite{Tsuge2006}. Chen \cite{Chen1997} proved a global existence theorem with large $L^\infty$ data having only non-negative initial velocity.  Recently, Huang-Li-Yuan \cite{Huang20171} proved the global existence of weak solutions of the isentropic Euler equations with spherical symmetry including the origin by $L^{\infty}$ compensated compactness framework and the vanishing viscosity method. If the blast wave moves outwards initially  and the initial densities and velocities decay to zero with certain rates near origin, then for any positive time, the densities and velocities tend to zero with the same rates near the origin. Chen et al.\cite{ChenP,ChenM} proved the global existence of finite-energy entropy solution by using $L^p$ compensated compactness framework. For the cylindrical symmetric flow, Chen-Glimm \cite{Chen3} considered the case excluding the origin and extended the analysis in \cite{Chen2} to $3\times3$ system. Chen-Wang \cite{ChenWang} applied the shock capturing schemes to the compressible Euler-Possion equations with geometrical structure in semiconductor devices. More interesting and relevent results can be found in \cite{Chen1997, Chen1998, Chen2003,Li,Liu,Tsuge2004a,Yang1995,Yang1996} and references therein.

 Almost all the results in previous work developed numerical scheme to construct approximate solutions and used lengthy estimates. In this paper, we apply the vanishing viscosity method together with the Smoller's invariant region theory \cite{Smoller} of quasilinear parabolic equations tackling with nonlinear singular terms to get uniform a priori estimates of approximate solutions. There are two technical difficulties to be overcome. One of the difficulties for the spherical or cylindrical symmetry flow is the singularity formed at the origin $x=0.$ There are complicated resonance between the characteristic mode and the geometric mode. To overcome this difficulty, we first employ the spacial scaling transformation $\rho=\bar{\rho}x^c,m=\bar{m}x^d,\tilde{m}=\hat{m}x^d$ as in Huang-Li-Yuan \cite{Huang20171}. Next, we introduce a new variable $\xi=\ln x$ and transform the original system into \eqref{geometriccc}. Then, viscous perturbations  $(\varepsilon \bar{\rho}_{\xi\xi}, \varepsilon \bar{m}_{\xi\xi},\varepsilon\hat{m}_{\xi\xi})$ are added in the new system for $(\bar{\rho},\bar{m},\hat{m})$. Finally, we get uniform estimates of approximate solutions independent of viscosity $\varepsilon$ using an enhanced maximum principle for parabolic system (See Lemma \ref{initial maximum}). We remark that Lemma \ref{initial maximum} is quite powerful when we estimate the parabolic systems with source terms.

 Another difficulty is from the third equation in \eqref{geometric}, due to the size of the system being $3\times3.$ There is an extra direction (tangential direction) we need to deal with. The extra field is linearly degenerate generating from the angular momentum. Different from spherical symmetric flow, extra estimate is needed on angular momentum $\tilde{m},$ which satisfies a parabolic equation. See \eqref{omega}. In fact, when $\tilde{m}=0,$ the equations could be reduced to the gas dynamics with cylindrically symmetric flow.  Motivated by the idea of quasi-decoupling method (See Appendix \ref{Appendix}) and framework developed by Chen \cite{Chen1992} and the numerical method furthered developed in Chen-Glimm \cite{Chen3} and Chen-Wang \cite{ChenWang}, as the case excluding the origin, we shall first obtain the strong convergence of approximate solutions $\rho^{\varepsilon},m^{\varepsilon}$ by using the estimates \eqref{solution3} and then we get strong convergence of the tangential momentum $\tilde{m}^{\varepsilon}.$ The quasi-decoupling technique established in \cite{Chen1992} is useful for studying the limiting behavior of approximate solutions to general hyperbolic conservation laws. As pointed in \cite{Chen1992}, initial oscillations could propagate along the corresponding third linearly degenerate field.  Our new observation is that if the blast wave initially moves outwards, that is, the normal velocity is non-negative and no restriction is needed on the sign of initial angular velocity, then both the density and velocities are continuous at the origin. Besides, if densities and velocities initially tend to zero at certain algebraic rate near the origin, then the density and velocity tend to zero at the same rate for any positive time.  \\

  We define the weak entropy solution and state our main result as follows.
\begin{Definition} \label{def3}
A measurable vector function $v(x, t)=(\rho(x,t), m(x,t),\tilde{m}(x,t))$ is called a global entropy solution of the Cauchy problem \eqref{geometricc2} provided that the following conditions hold:

(i) For any test function $\psi(x,t)=\varphi(x,t)A(x)\in C^1_0([0,\infty)\times[0,\infty)),$ $\varphi\in C^1_0([0,\infty)\times \R^+),$ the mass equation holds in the sense:
\begin{equation}\label{massweak} \int^{\infty}_0\int^{\infty}_0(\rho\psi_t+m\psi_x-\frac{1}{x}m\psi)dxdt+\int^{\infty}_0\rho_0(x)\psi(x,0)dx=0;
\end{equation}

(ii) For any test function $\psi_1(x,t)=\varphi_1(x,t)A(x)\in C^1_0([0,\infty)\times[0,\infty)),$ with $\varphi_1(x,t)\frac{\vec{x}}{x}\in C^1_0(\R^2\times \R^+),$ the normal momentum equation holds in the sense:
\begin{equation}\label{momentumweak} \int^{\infty}_0\int^{\infty}_0m\psi_{1t}+(\frac{m^2}{\rho}+p(\rho))\psi_{1x}-\frac{1}{x}\frac{m^2-\tilde{m}^2}{\rho}\psi_1dxdt+\int^{\infty}_0m_0(x)\psi_1(x,0)dx=0;
\end{equation}

(iii) For any test function $\psi_2(x,t)=\varphi_2(x,t)A(x)\in C^1_0([0,\infty)\times[0,\infty)),\varphi_2(x,t)\frac{(x_2,-x_1)}{x}\in C^1_0(\R^2\times \R^+),$ the angular momentum equations hold in the sense:
\begin{equation}\label{momentum2weak} \int^{\infty}_0\int^{\infty}_0\tilde{m}\psi_{2t}+\frac{m\tilde{m}}{\rho}\psi_{2x}-\frac{2m\tilde{m}}{x\rho}\psi_2dxdt+\int^{\infty}_0\tilde{m}_0(x)\psi_2(x,0)dx=0.
\end{equation}

For the weak entropy pair $(\eta, q)$, the inequality
\begin{equation*}
\eta(v)_t+q(v)_x-\nabla\eta(v)\cdot G(x,v)\leq0.
\end{equation*}
holds in the sense of distributions.
\end{Definition}
\begin{Remark}
 The entropy solutions obtained in our paper exactly satisfy the notions of the weak solutions of multidimensional rotating cylindrically symmetry flows in the sense of distributions. We shall show that weak solutions obtained  in the following Theorem \ref{include} are exactly the solutions of compressible rotating flow with cylindrical symmetry including the origin. It suffices to prove that the functions $(\rho(\vec{x},t),\vec{m}(\vec{x},t)),$ determined by the function $(\rho(x,t),m(x,t),\tilde{m}(x,t))$ of Theorem \ref{include}:
 \begin{equation}\label{structure}
 (\rho,\vec{m})(\vec{x},t)=\left(\rho(x,t),m(x,t)\frac{\vec{x}}{x}+\tilde{m}(x,t)\frac{(-x_2,x_1)}{x}\right),
 \end{equation}
 satisfy the standard definition of the weak solutions in which the origin is included. 
For any test function $\varphi\in C^1_0(\R^2\times \R^+)$ satisfying $\varphi(\vec{x},t)=\varphi(x,t),$ we set $\psi(\vec{x},t)=\varphi(\vec{x},t)A(x)\in C^1_0(\R^2\times[0,\infty)),$
 and obtain that
 \begin{equation}\label{mass}
 \begin{split}
 &\int^{\infty}_0\int_{\R^2}(\rho\varphi_t+\vec{m}\cdot\nabla\varphi)d\vec{x}dt+\int_{\R^2}\rho_0(\vec{x})\varphi(\vec{x},0)d\vec{x}\\
 &=\int^{\infty}_0\int^{\infty}_0(\rho\psi_t+m\psi_x-\frac{1}{x}m\psi)dxdt+\int^{\infty}_0\rho_0(x)\psi(x,0)dx=0.
 \end{split}
 \end{equation}
For any test function $\varphi_1\in C^1_0(\R^2\times\R^+),$ $\vec{\varphi}(\vec{x},t)=\varphi_1(x,t)\frac{\vec{x}}{x}\in C^1_0(\R^2\times \R^+),$
then we have
\begin{equation}\label{momentum}
 \begin{split}
 &\int^{\infty}_0\int_{\R^2}(\vec{m}\vec{\varphi}_{t}+\frac{\vec{m}\otimes \vec{m}}{\rho}\nabla\vec{\varphi}+p(\rho)div \vec{\varphi})d\vec{x}dt+\int_{\R^2}\vec{m}_0(\vec{x})\varphi(\vec{x},0)d\vec{x}\\
 =&2\pi\int^{\infty}_0\int^{\infty}_0xm\varphi_{1t}+x\frac{m^2}{\rho}\varphi_{1x}+p(\rho)(x\varphi_{1x}+\varphi_1)dxdt\\
&+2\pi\int^{\infty}_0 xm_0(x)\varphi_1(x,0)dx=0.
 \end{split}
 \end{equation}
 Setting $\psi_1(\vec{x},t)=\varphi_1(x,t)A(x)\in C^1_0([0,\infty)\times[0,\infty)),$ with $\varphi_1(x,t)\frac{\vec{x}}{x}\in C^1_0(\R^2\times \R^+),$ we obtain that
 \begin{equation}\label{momentum2}
  \int^{\infty}_0\int^{\infty}_0m\psi_{1t}+(\frac{m^2}{\rho}+p(\rho))\psi_{1x}-\frac{1}{x}\frac{m^2-\tilde{m}^2}{\rho}\psi_1dxdt+\int^{\infty}_0m_0(x)\psi_1(x,0)dx=0.
 \end{equation}
 For any test function $\varphi_2\in C^1_0(\R^2\times\R^+),$ $\vec{\varphi}(\vec{x},t)=\varphi_2(x,t)\frac{(x_2,-x_1)}{x}\in C^1_0(\R^2\times \R^+),$
then we have
\begin{equation}\label{momentum3}
 \begin{split}
 &\int^{\infty}_0\int_{\R^2}(\vec{m}\vec{\varphi}_{t}+\frac{\vec{m}\otimes \vec{m}}{\rho}\nabla\vec{\varphi}+p(\rho)div \vec{\varphi})d\vec{x}dt+\int_{\R^2}\vec{m}_0(\vec{x})\varphi(\vec{x},0)d\vec{x}\\
 =&2\pi\int^{\infty}_0\int^{\infty}_0x\tilde{m}\varphi_{2t}+\frac{m\tilde{m}}{\rho}(\varphi_2+x\varphi_{2x})-\frac{2m\tilde{m}}{\rho}\varphi_2dxdt\\
&+2\pi\int^{\infty}_0 x\tilde{m}_0(x)\varphi_2(x,0)dx=0.
 \end{split}
 \end{equation}
 Setting $\psi_2(\vec{x},t)=\varphi_2(x,t)A(x)\in C^1_0([0,\infty)\times[0,\infty)),$ with  $\varphi_2(x,t)\frac{(x_2,-x_1)}{x}\in C^1_0(\R^2\times \R^+),$ we obtain that
 \begin{equation}\label{momentum2}
  \int^{\infty}_0\int^{\infty}_0\tilde{m}\psi_{2t}+\frac{m\tilde{m}}{\rho}\psi_{2x}-\frac{2m\tilde{m}}{x\rho}\psi_2dxdt+\int^{\infty}_0\tilde{m}_0(x)\psi_2(x,0)dx=0.
 \end{equation}
\end{Remark}
\begin{Theorem}\text(Existence)\label{include}
Let $\gamma>1,$ and $d=(\theta+1)c>0, c=\frac{1}{\theta}.$ Assume that there exist two positive constants $M_1,M_2$ such that
\begin{equation}\label{ini2}
\begin{aligned}
&0\leq\rho_0(x),~\frac{m_0}{\rho_0}+\rho_0^{\theta} \leq M_1x,\,\frac{m_0}{\rho_0}-\rho_0^{\theta}\geq 0,\left|\frac{\tilde{m}_0}{\rho_0}\right|\leq M_2x, \text{ a.e. $x\in [0,+\infty)$}.\\
\end{aligned}
\end{equation}
Then there exists a global entropy solution of \eqref{geometric}-\eqref{sphericalinitial2} satisfying
\begin{equation}\label{solution}
0\leq\rho(x,t)\leq Cx^c,~0\leq m(x, t)\leq C\rho(x, t)x,~|\tilde{m}(x, t)|\leq C\rho(x, t)x,
\end{equation}
for \text{a.e.$(x,t)\in[0,+\infty)\times\R^{+}$}, where $C$ depends on $M_1,M_2,T.$
\end{Theorem}
\begin{Remark}
To the best of the author's knowledge, the entropy solution obtained above is the first result on the entropy solution of the Cauchy problem of isentropic gas dynamics with cylindrically symmetric rotating structure including the origin.
\end{Remark}
\begin{Remark}
The initial condition are allowed to be unbounded when $x$ is sufficiently large in Theorem \ref{include}.
\end{Remark}
The rest of the paper is organized as follows. In Section \ref{formula}, we describe our approach of constructing approximate solutions by adding artificial viscosity. Some necessary formulas for the system are also given in this section. In Section \ref{includeorigin}, the uniform upper bound estimate for the approximate solutions $(\rho^\eps, m^\eps,\tilde{m}^\eps)$ is proved and then the $H_{loc}^{-1}$ compactness of entropy pair, consistency and entropy inequality are obtained, and finally the proof of Theorem 1.1 is completed by applying $L^{\infty}$ compactness framework established in  \cite{Ding, Diperna, Lions2}. In Appendix \ref{Appendix}, we prove the strong convergence of angular momentum $\tilde{m}$ for completeness.

\section{Formulations}\label{formula}
First, we recall some basic notations for the system \eqref{cylindricalsym} with eigenvalues:
\begin{equation}\label{2.1}
\lambda_1=\frac{m}{\rho}-\theta\rho^\theta,\quad
\lambda_2=\frac{m}{\rho}+\theta\rho^\theta,\quad
\lambda_3=\frac{m}{\rho},\quad
\end{equation}
where $\theta=\frac{\gamma-1}{2},$\,and the corresponding right eigenvectors are
\begin{equation}\label{2.2}
r_1=\left[\begin{array}{cc}
1\\ \lambda_1 \\ -\frac{\tilde{m}}{\rho}
\end{array}
\right],\quad
r_2=\left[\begin{array}{cc}
1\\ \lambda_2\\-\frac{\tilde{m}}{\rho}
\end{array}
\right],\quad
r_3=\left[\begin{array}{cc}
0\\ 0\\1
\end{array}
\right].
\end{equation}
The Riemann invariants $(w, z,\omega)$ are given by
\begin{equation}\label{2.3}
w=\frac{m}{\rho}+\rho^\theta,\quad z=\frac{m}{\rho}-\rho^\theta,\omega=\frac{\tilde{m}}{\rho},
\end{equation}
satisfying
\begin{equation}\label{lineardegenerate}
\nabla w\cdot r_1=0,\nabla z\cdot r_2=0,\nabla \omega\cdot r_3=\frac{1}{\rho},
\end{equation}
where $\nabla=(\partial_\rho, \partial_m, \partial_{\tilde{m}})$ is the gradient with respect to the phase-space coordinates $v$.
 Moreover,
 \begin{equation}\label{lambdaa}
\nabla \lambda_1\cdot r_1=-\theta(1+\theta)\rho^{\theta-1},\nabla \lambda_2\cdot r_2=\theta(1+\theta)\rho^{\theta-1},\nabla \lambda_3\cdot r_3=0,
\end{equation}
 The third characteristic field $\lambda_3$ is linearly degenerate. From \eqref{lambdaa}, the first two characteristic field in \eqref{geometric} are genuinely nonlinear if $\gamma\in(1,3],$ the system is not genuinely nonlinear at vacuum $\rho=0$ if $\gamma>3.$ A pair of functions $(\eta, q):\R^3\mapsto\R^2$ is called an entropy-entropy flux pair of system \eqref{geometric} or \eqref{geometricc2} if it satisfies
\begin{equation}\label{2.4}
\nabla q(v)=\nabla\eta(v)\nabla\left[\begin{array}{ccc}
m\\ \frac{m^2}{\rho}+p(\rho) \\ \frac{m\tilde{m}}{\rho}
\end{array}
\right].
\end{equation}
When $$\eta\left|_{\rho=0, \frac{m}{\rho},\frac{\tilde{m}}{\rho}\text{ fixed }}=0,\right.$$ $\eta(\rho, m,\tilde{m})$ is defined to be weak entropy of system \eqref{geometric}.

The mechanical energy and mechanical energy flux, $(\eta^*,q^*),$
\begin{equation}\label{2.5}
\eta^*(\rho, m,\tilde{m})=\frac{m^2+\tilde{m}^2}{2\rho}+\frac{p_0\rho^\gamma}{\gamma-1},~~
q^*(\rho, m,\tilde{m})=\frac{m^3+m\tilde{m}^2}{2\rho^2}+\frac{\gamma p_0\rho^{\gamma-1}m}{\gamma-1},
\end{equation}
is a strictly convex weak entropy pair of system \eqref{geometric}.

Next, we will introduce a variant of invariant region theory for a decoupled parabolic system which can be used to deal with general nonlinear source terms. In the following, we extend the maximum principle to quasilinear parabolic systems and apply it to prove the global existence of isentropic compressible Euler equations with cylindrical symmetry. Readers could refer to \cite{Huang20171,Huang20172} for details.
\begin{Lemma}(Maximum principle)\label{initial maximum}
Let $p(x,t), q(x,t)$, $(x,t)\in[a,b]\times[0,T]$ be any bounded classical solutions of the following quasilinear parabolic system
\begin{eqnarray}\label{pq}
\left\{ \begin{split}
\displaystyle &p_t+\mu_1 p_x=
\eps p_{xx}+a_{11}p+a_{12}q+R_1,\\
\displaystyle &q_t+\mu_2 q_x=
\eps q_{xx}+a_{21}p+a_{22}q+R_2,
\end{split}
\right.
\end{eqnarray}
with initial-boundary data
\[\begin{aligned}
  p(x,0)&\leq 0,~q(x, 0)\geq0, \text{ for } ~x\in[a,b],\\
  p(a,t)& \leq 0,\,q(a,t)\geq0, \text{ for }~ t\in[0,T],\\
  p(b,t)& \leq 0,\,q(b,t)\geq0, \text{ for }~ t\in[0,T],\\
\end{aligned}\]
where $$\mu_{i}=\mu_i(x,t,p(x,t),q(x,t)),a_{ij}=a_{ij}(x,t,p(x,t),q(x,t)),$$and the source terms $$R_i=R_i(x,t,p(x,t),q(x,t),p_{x}(x,t),q_{x}(x,t)),i,j=1,2,\forall(x,t)\in[a,b]\times[0,T].$$ $\mu_{i},a_{ij}$are bounded with respect to $(x,t,p,q)\in[a,b]\times[0,T]\times K,$ where $K$ is an arbitrary compact subset in $\R^2.$  $a_{12},a_{21},R_{1},R_{2}$ are continuously differentiable with respect to $p,q.$\\
Assume the following conditions hold:\\
\begin{description}
  \item[(C1)]  When $p=0$ and $q\geq0,$ there is $a_{12}\leq0;$ When $q=0$ and $p\leq0,$  there is  $a_{21}\leq0;$\label{jj}
  \item[(C2)]  When $p=0$ and $q\geq0,$ there is $~R_1=R_1(x,t,0,q,\zeta,\eta)\leq0;$ When $q=0$ and $p\leq0,$ there is  $~R_2=R_2(x,t,p,0,\zeta,\eta)\geq0.$\\
\end{description}
\end{Lemma}
Then for any $(x, t)\in[a,b]\times[0,T],$
$$p(x,t)\leq 0, ~~q(x, t)\geq0.$$


\section{Proof of Theorem 1.1}\label{includeorigin}
\subsection{Uniform estimates}\label{upper}

Here, we consider \eqref{geometric} on a cylinder $(a,b)\times \R^{+},$ with $\R^{+}=[0,+\infty),a:=a(\eps)=-\frac{1}{\ln \eps},$  $\lim\limits_{\eps\rightarrow0}b(\eps)=\infty.$
We set the transformation to the equation \eqref{geometric}:
\begin{equation}\label{tra}
\rho=\bar{\rho}x^c,m=\bar{m}x^d,\tilde{m}=\hat{m}x^d.
\end{equation}
Motivated by \cite{Tsuge2004a}, we take $d=(\theta+1)c=c+1,c=\frac{1}{\theta}>0.$
Then we can transform \eqref{geometric} into following form:
\begin{eqnarray}\label{geometric222}
\left\{ \begin{split} \bar{\rho}_t+\bar{m}_xx=&-(d+1)\bar{m},\,\\ \bar{m}_t+\left(\frac{\bar{m}^2}{\bar{\rho}}+p(\bar{\rho})\right)_xx=&-(2d-c+1)\frac{\bar{m}^2}{\bar{\rho}}+\frac{\hat{m}^2}{\bar{\rho}}-(2d-c)p(\bar{\rho}),\\
\hat{m}_t+\left(\frac{\bar{m}\hat{m}}{\bar{\rho}}\right)_xx=&-(2d-c+2)\frac{\bar{m}\hat{m}}{\bar{\rho}},
~x\in [a(\eps),b(\eps)],t>0.
\end{split}
\right.
\end{eqnarray}
We introduce a new variable $\xi=\ln x.$
Then \eqref{geometric222} can be transformed into the following form:
\begin{eqnarray}\label{geometriccc}
\left\{ \begin{split} \bar{\rho}_t+\bar{m}_\xi=&-(d+1)\bar{m},\,\\ \bar{m}_t+\left(\frac{\bar{m}^2}{\bar{\rho}}+p(\bar{\rho})\right)_\xi=&-(2d-c+1)\frac{\bar{m}^2}{\bar{\rho}}+\frac{\hat{m}^2}{\bar{\rho}}-(2d-c)p(\bar{\rho}),\\
\hat{m}_t+\left(\frac{\bar{m}\hat{m}}{\bar{\rho}}\right)_\xi=&-(2d-c+2)\frac{\bar{m}\hat{m}}{\bar{\rho}},
~x\in [a(\eps),b(\eps)],t>0.
\end{split}
\right.
\end{eqnarray}
We approximate \eqref{geometriccc} by adding the artificial viscosity as follows:
\begin{eqnarray}\label{geometric3}
\left\{ \begin{split} \bar{\rho}_t+\bar{m}_\xi=&-(d+1)\bar{m}+\eps\bar{\rho}_{\xi\xi},\,\\ \bar{m}_t+\left(\frac{\bar{m}^2}{\bar{\rho}}+p(\bar{\rho})\right)_\xi=&-(2d-c+1)\frac{\bar{m}^2}{\bar{\rho}}+\frac{\hat{m}^2}{\bar{\rho}}-(2d-c)p(\bar{\rho})+\eps \bar{m}_{\xi\xi},\\
\hat{m}_t+\left(\frac{\bar{m}\hat{m}}{\bar{\rho}}\right)_\xi=&-(2d-c+2)\frac{\bar{m}\hat{m}}{\bar{\rho}}+\eps \hat{m}_{\xi\xi},\\
&x\in [a(\eps),b(\eps)],t>0.
\end{split}
\right.
\end{eqnarray}
The initial-boundary value conditions for \eqref{geometric3} are given by£º
\begin{equation}\label{geometricini-vis22}
\begin{aligned}
&(\bar{\rho}, \bar{m},\hat{m})|_{t=0}=(\bar{\rho}_0^\eps(x),\bar{m}_0^\eps(x),\hat{m}_0^\eps(x))\\
&=(\bar{\rho}_0(x)+\eps^{\frac{2}{\theta}}, (\frac{\bar{m}_0(x)}{\bar{\rho}_0(x)}+\eps)(\bar{\rho}_0(x)+\eps^{\frac{2}{\theta}})\chi,(\frac{\hat{m}_0(x)}{\bar{\rho}_0(x)}+\eps)(\bar{\rho}_0(x)+\eps^{\frac{2}{\theta}})\chi)\ast j^\eps,\\
& x\in [a(\eps),b(\eps)];\\
&(\bar{\rho},\bar{m},\hat{m})|_{x=a(\eps)}=(\bar{\rho}_0^\eps (a(\eps)), \bar{m}_0^\eps (a(\eps)),\hat{m}_0^\eps (a(\eps)))=(\bar{\rho}_0^\eps (a(\eps)),0,0),\\
&(\bar{\rho},\bar{m},\hat{m})|_{x=b(\eps)}=(\bar{\rho}_0^\eps (b(\eps)), \bar{m}_0^\eps (b(\eps)),\hat{m}_0^\eps (b(\eps))), t>0,
\end{aligned}
\end{equation}
where $j^\eps$ is the standard mollifier, $\chi=\chi_{[2a(\eps),b(\eps)]}$ is the characteristic function  and the parameter $\eps>0$ is small.

Next, we will derive the uniform bound of the approximate solution $\bar{\rho},$$\bar{m},$$\hat{m}$ by the maximum principle, i.e., Lemma \ref{initial maximum} and then derive the $L^\infty$ bound for $\rho,m,\tilde{m}.$
Using Riemann invariants's definition, we obtain that
$$w=\frac{m}{\rho}+\rho^\theta=\left(\frac{\bar{m}}{\bar{\rho}}+\bar{\rho}^\theta\right)x:=\bar{w}x,z=\frac{m}{\rho}-\rho^\theta=\left(\frac{\bar{m}}{\bar{\rho}}-\bar{\rho}^\theta\right)x:=\bar{z}x,$$
$$\omega=\frac{\tilde{m}}{\rho}=\frac{\hat{m}}{\bar{\rho}}x:=\hat{\omega}x.$$
Similarly, we have $$\lambda_1=\frac{m}{\rho}-\theta\rho^\theta=\left(\frac{\bar{m}}{\bar{\rho}}-\theta\bar{\rho}^\theta\right)x:=\bar{\lambda}_1x,\lambda_2=\frac{m}{\rho}+\theta\rho^\theta=\left(\frac{\bar{m}}{\bar{\rho}}+\theta\bar{\rho}^\theta\right)x:=\bar{\lambda}_2x,$$
$$\lambda_3=\frac{m}{\rho}=\frac{\bar{m}}{\bar{\rho}}x:=\bar{\lambda}_3x.$$

Then we can approximate \eqref{geometric} by adding the same artificial viscosity:
\begin{eqnarray}\label{geometricvis2}
\left\{ \begin{split}
\displaystyle \rho_t+m_x=&-\frac{1}{x}m+\eps\left[\rho_{xx}x^{2}+(d-3c)\rho_xx+c^2\rho \right].\\
 m_t+\left(\frac{m^2}{\rho}+p(\rho)\right)_x=&-\frac{1}{x}\frac{m^{2}-\tilde{m}^2}{\rho}+\eps\left[m_{xx}x^{2}-(c+d)m_xx+d^2m\right]\\
 \tilde{m}_t+\left(\frac{m\tilde{m}}{\rho}\right)_x=&-\frac{2}{x}\frac{m\tilde{m}}{\rho}+\eps\left[\tilde{m}_{xx}x^2-(c+d)\tilde{m}_xx+d^2\tilde{m}\right],x\in[a(\eps),b(\eps)],t>0,
\end{split}
\right.
\end{eqnarray}

It is noted that system \eqref{geometric3} is equivalent to system \eqref{geometricvis2}. Therefore it suffice to prove the existence of the approximate solutions in \eqref{geometricvis2}.\\

Now, we transform \eqref{geometric3} into the following form:
\begin{eqnarray}\label{wwzz222}
\left\{ \begin{split}
\bar{w}_t+\bar{\lambda}_2 \bar{w}_\xi=&\eps \bar{w}_{\xi\xi}+2\eps \frac{\bar{\rho}_\xi}{\bar{\rho}}\bar{w}_\xi-\eps\theta(\theta+1)\bar{\rho}^{\theta-2}\bar{\rho}_\xi^2\\
&-\frac{\bar{m}^2}{\bar{\rho}^2}-\theta(d+1)\bar{\rho}^{\theta-1}\bar{m}-(2d-c)\frac{\theta^2}{\gamma}\bar{\rho}^{2\theta}+\frac{\hat{m}^2}{\bar{\rho}^2},\\
 \bar{z}_t+\bar{\lambda}_1  \bar{z}_\xi=&\eps\bar{z}_{\xi\xi}+2\eps \frac{\bar{\rho}_\xi}{\bar{\rho}}\bar{z}_\xi+\eps\theta(\theta+1)\bar{\rho}^{\theta-2}\bar{\rho}_\xi^2\\
&-\frac{\bar{m}^2}{\bar{\rho}^2}+\theta(d+1)\bar{\rho}^{\theta-1}\bar{m}-(2d-c)\frac{\theta^2}{\gamma}\bar{\rho}^{2\theta}+\frac{\hat{m}^2}{\bar{\rho}^2},\\
\hat{\omega}_t+\bar{\lambda}_3\hat{\omega}_\xi=&\eps\hat{\omega}_{\xi\xi}+2\eps \frac{\bar{\rho}_\xi}{\bar{\rho}}\hat{\omega}_\xi-2\frac{\bar{m}\hat{m}}{\bar{\rho}^2}.\\
\end{split}
\right.
\end{eqnarray}

We first consider the $L^\infty$ estimate of $\hat{\omega}$ satisfying  the third equation in \eqref{wwzz222},

\begin{equation}\label{omega}
\hat{\omega}_t+\left(\bar{\lambda}_3-2\eps \frac{\bar{\rho}_\xi}{\bar{\rho}}\right)\hat{\omega}_\xi=\eps\hat{\omega}_{\xi\xi}-2\frac{\bar{m}}{\bar{\rho}}\hat{\omega}.
\end{equation}
 From the initial assumption \eqref{ini22} on $\hat{\omega},$ there hold
 $$\Big|\frac{\tilde{m}_0}{\rho_0}\Big|\leq M_2,\Big|\frac{\tilde{m}}{\rho}(a(\eps),t)\Big|\leq M_2 a(\eps),\Big|\frac{\tilde{m}}{\rho}(b(\eps),t)\Big|\leq M_2 b(\eps).$$
 Then, we have
  $$|\hat{\omega}_0|\leq M_2,|\hat{\omega}(a(\eps),t) |\leq M_2,|\hat{\omega}(b(\eps),t)|\leq M_2.$$
Then using the standard maximum principle for \eqref{omega}, we have
\begin{equation}\label{omegaestimate}
||\hat{\omega}||_{L^\infty}\leq C,
\end{equation}
where $C$ is dependent on initial data $||\rho_0,\tilde{m}_0||_{L^\infty},$ which is independent of $\eps.$ We will prove the $L^{\infty}$ bound of $\rho,m$ by Lemma \ref{initial maximum} and by the positiveness of $\frac{\bar{m}}{\bar{\rho}}$.
We define the control functions $(\phi,\psi)=(M_3+Ct+2\eps,0),$  and the modified Riemann invariants $(\hat{w},\hat{z})$ as
\begin{equation}\label{r}
\hat{w}=\bar{w}-\phi, ~~\hat{z}=\bar{z}+\psi.
\end{equation}
The system \eqref{wwzz222} becomes
\begin{eqnarray}\label{wwzz33}
\displaystyle\left\{ \begin{split}
\hat{w}_t+\bar{\lambda}_2 \hat{w}_\xi=&\eps \hat{w}_{\xi\xi}+2\eps\frac{\bar{\rho}_\xi}{\bar{\rho}}\hat{w}_\xi-\eps\theta(\theta+1)\bar{\rho}^{\theta-2}\bar{\rho}_{\xi}^2\\
+&\left[-\frac{\bar{m}^2}{\bar{\rho}^2}-\theta(d+1)\bar{\rho}^{\theta}\hat{z}-\theta(d+1)\bar{\rho}^{2\theta}-\theta^2c\bar{\rho}^{2\theta}+\hat{\omega}^2-C\right],\\
\hat{z}_t+\bar{\lambda}_1\hat{z}_\xi=&\eps \hat{z}_{\xi\xi}+2\eps\frac{\bar{\rho}_\xi}{\bar{\rho}}\hat{z}_\xi+\eps\theta(\theta+1)\bar{\rho}^{\theta-2}\bar{\rho}_{\xi}^2\\
+&\left[-(\frac{\bar{w}+\hat{z}}{2})^2+\theta(d+1)\frac{\bar{w}^2-\hat{z}^2}{4}-\theta^2 c\left(\frac{\bar{w}-\hat{z}}{2}\right)^2+\hat{\omega}^2\right].
\end{split}
\right.
\end{eqnarray}
Then the above system \eqref{wwzz33} can be written as
\begin{eqnarray}\label{rst}
\displaystyle\left\{ \begin{split} &\hat{w}_t+\left(\bar{\lambda}_2-2\eps\frac{\bar{\rho}_\xi}{\bar{\rho}} \right)\hat{w}_\xi
=\eps\hat{w}_{\xi\xi}+a_{11}\hat{w}
+a_{12}\hat{z}+R_1,\\
&\hat{z}_t+\left(\bar{\lambda}_1-2\eps\frac{\bar{\rho}_\xi}{\bar{\rho}} \right) \hat{z}_\xi
=\eps\hat{z}_{\xi\xi}+a_{21}\hat{w}
+a_{22}\hat{z}+R_2,
\end{split}
\right.
\end{eqnarray}
where
\begin{equation*}
\begin{split}
&a_{11}=0,\,a_{12}=-\theta(d+1)\bar{\rho}^\theta\leq0,\,\\
&a_{21}=0,\,a_{22}=\frac{1}{4}\left[-1-\theta(d+1)-\theta^2c\right]\hat{z}+\frac{1}{2}\left(\theta^2c-1\right)\bar{w},\\
&R_1=-\frac{\bar{m}^2}{\bar{\rho}^2}-\theta(d+1)\bar{\rho}^{2\theta}-\theta^2c\bar{\rho}^{2\theta}-\eps\theta(\theta+1)\bar{\rho}^{\theta-2}\bar{\rho}_{\xi}^2+(\hat{\omega}^2-C)\leq0,\\
\,&R_2=\frac{1}{4}\left[-1+\theta(d+1)-\theta^2c\right]\bar{w}^2+\eps\theta(\theta+1)\bar{\rho}^{\theta-2}\bar{\rho}_{\xi}^2+\hat{\omega}^2\geq\frac{1}{4} \theta\bar{w}^2\geq0.
\end{split}
\end{equation*}
By the initial and boundary data,
\begin{equation}\label{initial}
\begin{split}
&\frac{m_0}{\rho_0}+\rho_0^{\theta}\leq (M_3+2\eps)x,\frac{m_0}{\rho_0}-\rho^{\theta}_0\geq0, x\in[a(\eps),b(\eps))],\\
&\left(\frac{m}{\rho}+\rho^{\theta}\right)|_{x=a(\eps)}\leq (M_3+2\eps) a(\eps),
\left(\frac{m}{\rho}-\rho^{\theta}\right)|_{x=a(\eps)}\geq0,\\
&\left(\frac{m}{\rho}+\rho^{\theta}\right)|_{x=b(\eps)}\leq (M_3+2\eps) b(\eps),
\left(\frac{m}{\rho}-\rho^{\theta}\right)|_{x=b(\eps)}\geq0,
\end{split}
\end{equation}
we obtain that
$$\hat{w}(\xi,0)=\bar{w}(\xi,0)-M_3-2\eps-Ct\leq 0,\,\hat{z}(\xi, 0)=\bar{z}(\xi, 0)\geq0,\ln(a(\eps))\leq\xi\leq\ln(b(\eps)),$$
$$\bar{w}(\ln(a(\eps)),t)\leq (M_3+2\eps) ,\,\bar{z}(\ln(a(\eps)),t)=0,\,\text{for } t>0,\,$$
$$\bar{w}(\ln(b(\eps)),t)\leq (M_3+2\eps) ,\,\bar{z}(\ln(b(\eps)),t)=0,\,\text{for } t>0.\,$$
By Lemma \ref{initial maximum}, $$\bar{w}(\xi,t)\leq C(T), \bar{z}(\xi,t)\geq 0.$$
Hence, we obtain $$0\leq\bar{\rho}(x,t)\leq C(T), 0\leq \bar{m}(x,t)\leq C\bar{\rho}(x,t),|\hat{m}(x,t)|\leq C\rho(x,t),$$ i.e.,  \text{ for  a.e.$(x,t)\in[a(\eps),b(\eps)]\times\R^{+}$},
\begin{equation}\label{solution3}
0\leq\rho(x,t)\leq Cx^c,~0\leq m(x, t)\leq C\rho(x, t)x, |\tilde{m}(x,t)|\leq C\rho(x,t)x,
\end{equation}
where $C$ depends on $M_3,T.$
 This completes the proof of the following theorem.
 \begin{Theorem}\text($L^{\infty}$ estimate for Cylindrically Symmetry Problem)\label{main3} Let $\gamma>1.$ For any positive constants $c$ and $d$ satisfying $d=(\theta+1)c>0,c=\frac{1}{\theta}.$ Assume that there exists a positive constant $M_3$ such that the initial and boundary data satisfy
\begin{equation}\label{ini22}
\begin{aligned}
&\rho_0(x)\geq \eps^{\frac{2}{\theta}}x^c,~\frac{m_0}{\rho_0}+\rho_0^{\theta} \leq (M_3+2\eps) x,\,\frac{m_0}{\rho_0}-\rho_0^{\theta}\geq 0,\left|\frac{\tilde{m}_0}{\rho_0}\right|\leq M_2,\text{ a.e. $x\in [a(\eps),b(\eps)]$},\\
&\left(\frac{m}{\rho}+\rho^{\theta}\right)|_{x=a(\eps)}\leq (M_3+2\eps) a(\eps),\left(\frac{m}{\rho}-\rho^{\theta}\right)|_{x=a(\eps)} \geq0,\\
&\left(\frac{m}{\rho}+\rho^{\theta}\right)|_{x=b(\eps)}\leq (M_3+2\eps) b(\eps),\left(\frac{m}{\rho}-\rho^{\theta}\right)|_{x=b(\eps)} \geq0,\\
&\left|\frac{\tilde{m}}{\rho}(a(\eps),t)\right|\leq M_2a(\eps),\left|\frac{\tilde{m}}{\rho}(b(\eps),t)\right|\leq M_2b(\eps), \,t>0.\\
\end{aligned}
\end{equation}
Then for the solution of \eqref{geometricvis2}, \eqref{initial}, the following estimates hold
\begin{equation}\label{solution2}
0\leq\rho(x,t)\leq Cx^c,~0\leq m(x, t)\leq C\rho(x, t)x, |\tilde{m}(x,t)|\leq C\rho(x,t)x,
\end{equation}
for a.e.,
$(x,t)\in[a(\eps),b(\eps)]\times\R^{+},$ where $C$ depends on $M_2,M_3,T.$
\end{Theorem}
\subsection{Lower bound estimate}\label{low2}

From \eqref{solution2} , we can obtain the lower bound of the density and the global existence of approximate solutions, whose proof is similar to that in \cite{Huang20171,Huang20172}.\\

\begin{Theorem}\label{app}
For any time $T>0,$ there exist positive constants $C$ and $\eps_0$  such that for $0<\eps<\eps_0$, the initial-boundary value problem
 \eqref{geometricvis2}, \eqref{initial} has following $L^{\infty}$ estimates
\begin{equation*}
e^{-C(\eps, T)}x^c\leq \rho^\eps(x, t)\leq Cx^c, ~~0\leq m^\eps(x, t)\leq C\rho^\eps(x, t)x, |\tilde{m}^{\eps}(x,t)|\leq C\rho^{\eps}(x,t)x,
\end{equation*}
where $C$ is independent of $x.$
\end{Theorem}

\subsection{$H^{-1}_{loc}$ compactness of the entropy pair}\label{entropy}
 Denote $\Pi_T=(a(\eps),b(\eps))\times [0, T],$ for any $T\in(0, \infty).$  Let $K\subset\Pi_T$ be any compact set, and choose $\varphi\in C_c^\infty(\Pi_T)$ such that $\varphi|_{K}=1,$ and $0\leq\varphi\leq1.$
Multiplying \eqref{geometric3} by  $\nabla\eta^* \varphi$ with $\eta^*$ the mechanical entropy,
we obtain
\begin{equation}\label{entropy1}
\begin{split}
&\eps\int\int_{\Pi_T}\varphi(\rho_x, m_x,\tilde{m}_x)\nabla^2\eta^*(\rho_x, m_x,\tilde{m}_x)^\top x^{2}dxdt\\
=&\int\int_{\Pi_T}\left[\eta^*_\rho\left(-\frac{1}{x}m+\eps(d-3c)\rho_xx+\eps c^2\rho\right)\right.\\
&+\left.\eta^*_m\left(-\frac{1}{x}\frac{m^2}{\rho}-\eps(c+d) m_xx+\eps d(c+1)m \right)\right.\\
&+\left.\eta^*_{\tilde{m}}\left(-\frac{2}{x}\frac{m\tilde{m}}{\rho}-\eps(c+d) \tilde{m}_xx+\eps d(c+1)\tilde{m} \right)\right]\varphi\\
&+\eta^*\varphi_t+q^*\varphi_x+\eps\eta^*(x^{2}\varphi)_{xx}dxdt.\\
\end{split}
\end{equation}
Direct calculation shows that
\[
(\rho_x, m_x,\tilde{m}_x)\nabla^2\eta^*(\rho_x, m_x,\tilde{m}_x)^\top=p_0\gamma\rho^{\gamma-2}\rho_x^2
+\rho (u_x^2+\tilde{u}_x^2).
\]
Note that
\begin{equation*}
|\eta^*_m\eps m_xx(c+d)|\leq\frac{\eps p_0\gamma}{8}\rho^{\gamma-2}\rho_x^2x^2+C_{\eps}\rho^{2-\gamma}u^4+\frac{\eps \rho u_x^2}{2}x^2+C_\eps\rho u^2,\\
\end{equation*}
$$|\eta^*_\rho\eps \rho_x(d-3c)x|\leq\frac{\eps p_0\gamma}{4}\rho^{\gamma-2}\rho_x^2x^{2}+C_{\eps}\rho^{2-\gamma}u^4+C_\eps\rho^{2(\gamma-1)},$$
\begin{equation*}
|\eta^*_{\tilde{m}}\eps \tilde{m}_xx(c+d)|\leq\frac{\eps p_0\gamma}{8}\rho^{\gamma-2}\rho_x^2x^2+C_{\eps}\rho^{2-\gamma}u^4+\frac{\eps \rho \tilde{u}_x^2}{2}x^2+C_\eps\rho \tilde{u}^2,\\
\end{equation*}
we get
\begin{equation*}
\begin{split}
&\frac{\eps}{2}\int\int_{\Pi_T}\varphi(\rho_x, m_x,\tilde{m}_x)\nabla^2\eta^*(\rho_x, m_x,\tilde{m}_x)^\top x^{2} dxdt\\
\leq&\int\int_{\Pi_T}\left(C_{\eps}\rho^{2-\gamma}u^4+C_{\eps}\rho u^2+C_{\eps}\rho^{2(\gamma-1)}\right)\varphi dxdt\\
&+\int\int_{\Pi_T}\eta^*\varphi_t+q^*\varphi_x+\eps\eta^*(x^{2}\varphi)_{xx}+(\frac{m^{2}+\tilde{m}^2}{2\rho^{2}}
-\frac{p_0\gamma}{\gamma-1}\rho^{\gamma-1})\frac{1}{x}m\varphi-\frac{1}{x}\frac{m^3}{\rho^2}\varphi dxdt\\
\leq& C(\varphi),\\
\end{split}
\end{equation*}
where the constant $C(\varphi)$ depends on $\varphi$. Thus we have arrived that
$$\eps(\rho_x, m_x,\tilde{m}_x)\nabla^2\eta^*(\rho_x, m_x,\tilde{m}_x)^\top\in L^1_{loc}(\Pi_T),$$
i.e.,
\begin{equation}\label{locestimate}
 \eps\rho^{\gamma-2}\rho_x^2+\eps\rho(u_x^2+\tilde{u}_x^2)\in L^1_{loc}(\Pi_T).
\end{equation}
For any weak entropy-entropy flux pairs independent of $\tilde{m}$, then we have
\begin{equation}\label{4.3}
\begin{split}
\eta_t+q_x=&\eps\eta_{xx}x^2-\eps(\rho_x, m_x)\nabla^2\eta(\rho_x, m_x)^\top x^2\\
+&\left(\eta_\rho(-\frac{1}{x}m+\eps c^2\rho)+\eta_m(-\frac{1}{x}\frac{m^{2}}{\rho}+\eps d(c+1)\tilde{m})\right)\\
+&\eps\left(\eta_\rho\rho_x(d-3c)-\eta_mm_x(c+d)\right)=:\sum_{i=1}^4J_i.
\end{split}
\end{equation}
From \eqref{locestimate}, it is obvious that
$J_1$ is compact in $H^{-1}_{loc}(\Pi_T).$ The proof is similar with the argument in \cite{Lu} and will not be reproduced here.
For $\gamma>2,\rho^{\gamma-2}\rho^2_x$ is degenerate near the vacuum. Here, we assume that $1<\gamma\leq2.$
For any weak entropy, the Hessian matrix $\nabla^2\eta$ can be controlled by  $\nabla^2\eta^*$ , refer to Lions et al. \cite{Lions2},
\begin{equation}\label{4.4}
(\rho_x, m_x)\nabla^2\eta(\rho_x, m_x)^\top\leq(\rho_x, m_x,\tilde{m}_x)\nabla^2\eta^*(\rho_x, m_x,\tilde{m}_x)^\top,
\end{equation}
therefore $J_2$ is bounded in $L^1_{loc}(\Pi_T)$ and by the Sobolev embedding theorem, this term is compact in
$W_{loc}^{-1, \alpha}(\Pi_T)$ for some $1<\alpha<2$.
For $J_3$,  $$|J_3|=\left|\eta_\rho(-\frac{1}{x}m+\eps c^2\rho)+\eta_m(-\frac{1}{x}\frac{m^{2}}{\rho}+\eps d(c+1)\tilde{m})\right|\leq C,$$
we obtain that $J_3$ is bounded in $L^1_{loc}(\Pi_T).$
For the last term $J_4$, we get
$$|J_4|\leq C\eps\left(\rho^{\gamma/2-1}|\rho_x|+\rho^{\frac{1}{2}}(|u_x|+|\tilde{u}_x|)\right)x.$$ Using  \eqref{locestimate}, we have $J_4$ is compact in $H^{-1}_{loc}(\Pi_T).$
Therefore,
\begin{equation}\label{H-1}
\eta_t+q_x \text{ is compact in } W^{-1, \alpha}_{loc}(\Pi_T) \text{ for some } 1<\alpha<2.
\end{equation}
It is easy to see that
\begin{equation}\label{H-11}
\eta_t+q_x \text{ bounded in } W^{-1, \infty}_{loc}(\Pi_T).
\end{equation}
Next, we shall utilize Murat's Lemma to conclude our result.
\begin{Lemma}{(Murat \cite{Murat})}\label{murat}
	Let $\Omega\subseteq\R^n$ be an open set, then
	\begin{equation*}
		(\text{compact set of } W^{-1, q}_{loc}(\Omega))\cap(\text{bounded set of } W^{-1, r}_{loc}(\Omega))\\
		\subset(\text{compact set of } H^{-1}_{loc}(\Omega)),
	\end{equation*}
	where $1<q\leq 2<r.$
\end{Lemma}
Combining Lemma \ref{murat}, \eqref{H-1} and \eqref{H-11}, we obtain that
\begin{equation}\label{compact}
\eta_t+q_x~ \text{is compact in}~ H^{-1}_{loc}(\Pi_T)
\end{equation}
holds for all weak entropy-entropy flux pairs, independent of $\tilde{m}.$
\begin{Remark}
 \eqref{compact} still holds for the case that $\gamma>2$ by a similar discussion as in  \cite{Wang,Lu}. We assume $1<\gamma\leq2$ for simplicity.
\end{Remark}
\subsection{Strong convergence and consistency}\label{Strong}
Combining the compactness in \eqref{compact} and the compactness framework established in
 \cite{Ding, Diperna, Lions2,Chen3},  we obtain that strong convergency of the approximate solutions. There exists
 a subsequence of $(\rho^\eps,m^\eps,\tilde{m}^\eps)$ (still adopt the same notation by $(\rho^\eps,m^\eps,\tilde{m}^\eps)$) such that
 \begin{equation}
 \label{4.16}
 (\rho^\eps, m^\eps,\tilde{m}^\eps)\to(\rho, m,\tilde{m}) ~~~
 \text{ in } L^p_{loc}(\Pi_T), ~~p\geq1. \end{equation}
 Here, we first apply the standard $2\times2$ compensated compactness framework established by \cite{Diperna,Ding,Ding1,Lions1} to prove the strong convergence of $(\rho^\eps,m^\eps).$ Next, we apply the quasi-decoupling method and framework for the limiting behavior of approximate solutions in \cite{Chen1992}, we can prove the strong convergence of $\tilde{m}^{\eps}.$ The prove is similar as in \cite{Chen3}, see Appendix \ref{Appendix} for the proof.

Furthermore, it is shown that $\rho=O(1)x^c$ and $m=O(1)x^{d},\tilde{m}=O(1)x^{d}$ so that the terms on right hand sides of the equation \eqref{geometric}, that is, $\frac{m}{x}$ and $\frac{m^2-\tilde{m}^2}{x\rho},$ and $\frac{2m\tilde{m}}{x\rho}$  are integrable near the origin with respect to the spacial variable $x$.  As in \cite{Chen2, Chen2003},  we can show that $(\rho, m,\tilde{m})$ is an entropy  solution to the problem \eqref{geometric}.  Therefore, this completes the proof of Theorem \ref{include}.

\appendix
\section{Strong convergence of $\tilde{m}^{\eps}$}\label{Appendix}
By the compensated compactness framework \cite{Diperna,Ding,Ding1,Lions1}, we can prove that there exists a subsequence of $(\rho^{\eps}(x,t),m^{\eps}(x,t))\rightarrow(\rho(x,t),m(x,t)),$ a.e. in $L^{\infty}.$ The pair of functions $(\rho(x,t),m(x,t))$ satisfies \eqref{solution2} and
$$(x\rho(x,t))_t+(xm(x,t))_x=0$$
in the sense of distributions with respect to the generalized test function space
$C^1_0(\Pi_T-\{\rho(x,t)=0\})\cap BV(\Pi_T).$
To prove the strong convergence of $\tilde{m}^{\eps}(x,t),$ we adopt the method of quasi-decoupling method introduced in \cite{Chen1992} and present here for convenience.

Consider a family of Radon measures $\{\mu_{x,t}:\Pi_T\rightarrow Prob(\R)\}$
with supp $\mu_{x,t}\subseteq\R$ and a family of initial Radon measures ${\mu_{x,0}}$
satisfying
\begin{eqnarray}\label{measure}
\left\{ \begin{array}{ll}
\displaystyle (<\mu_{x,t},a(x\omega-k)>P(x,t))_t+(<\mu_{x,t},a(x\omega-k)>\lambda(x,t)P(x,t))_x\leq0,\\
\displaystyle \mu_{x,t}|_{t=0}=\mu_{x,0},\\
\end{array}
\right.
\end{eqnarray}
in the generalized sense of distributions in $\Pi_T$ for any convex functions $a\in C(\R)$ and $k\in \R.$ Here, $(\lambda(x,t)P(x,t),P(x,t)),$ where $P(x,t)\geq0,$ is the divergence-free vector field in $(x,t)$ space:
$$div_{x,t}(\lambda(x,t)P(x,t),P(x,t))=0.$$
\begin{Lemma}{(Quasi-decoupling method\cite{Chen1992})}\label{quasidecouple}
	Suppose that a family of probability measures ${\mu_{x,t}:\Pi_T\rightarrow Prob(\R)}$ with compact support in $\R$ satisfy \eqref{measure} with a corresponding family of initial Dirac measures ${\mu_{x,0}=\delta_{\omega_0(x)}}$ for any convex function $a\in C^2(\R)$ and $k\in \R.$ Then,
\begin{equation}\nonumber
\mu_{x,t}(\omega)=\delta_{\omega_0(T^{-1}_0(p(x,t)))}(\omega), a.e., in \{(x,t)\in \Pi_T:P(x,t)>0\},
\end{equation}
where $p:\Pi_T\rightarrow \R^1$ is determined by the quasi-decoupling transformation,
\begin{equation}\nonumber
\frac{\partial p(x,t)}{\partial x}=P(x,t), \frac{\partial p(x,t)}{\partial t}=-\lambda(x,t)P(x,t),p(x,0)=\int^{x}_0P(\xi,0)d\xi,
\end{equation}
and $T^{-1}_0(p(x,0))=x, T^{-1}_0(p)$ is the inverse map of $T_t(x)|_{t=0}:=p(x,0).$
\end{Lemma}

For the approximate solution constructed in Theorem \ref{app}, $v^{\eps}(x,t)=(\rho^{\eps}(x,t),m^{\eps}(x,t)=\rho^{\eps}u^{\eps}(x,t),\tilde{m}^{\eps}(x,t)=\rho^{\eps}\omega^{\eps}(x,t))$ satisfy
$$(x\rho^{\eps}a(x \omega^{\eps}-k))_t+(xm^{\eps}a(x\omega^{\eps}-k))_x\leq o_{\eps}(1)\rightarrow0,\eps\rightarrow 0$$
in the sense of distributions for any convex function $a\in C^2$ by using the estimates \eqref{locestimate}.
 For the Young measure-valued solutions $\nu_{x,t}$ determined by the appropriate solutions $(\rho^{\eps},m^{\eps},\tilde{m}^{\eps}),$ we define Radon measures $\mu_{x,t}$ as the projection measures of $\nu_{x,t}$ on the state space $\omega\in \R:$
 $$<\mu_{x,t},g(\omega;x,t)>\equiv<\nu_{x,t},g(\omega;x,t)>$$
  for any $g(\omega;x,t)\in C^0.$
  Then, we obtain that $\mu_{x,t}$ satisfy
  $$(x\rho(x,t)<\mu_{x,t},a(x\omega-k)>)_t+(xm(x,t)<\mu_{x,t},a(x\omega-k)>)_x\leq0,$$
   for any convex function $a\in C^2(\R)$ and $k\in \R,$ where the function $\rho(x,t)$ and $m(x,t)$ are the limit functions of the approximate solutions $\rho^{\eps},m^{\eps}.$
   We note that
   $$(x\rho(x,t))_t+(xm(x,t))_x=0$$
   satisfies in the sense of distributions and $x\rho(x,t)\geq0,$
   Furthermore, $$\omega^{\eps}_0(x)\rightarrow \omega_0(x), a.e.,$$
   this implies $$\mu_{x,0}(\omega)=\delta_{\omega_0(x)}(\omega).$$
Using Lemma \ref{quasidecouple}, we deduce that
 $$\mu_{x,t}(\omega)=\delta_{\omega(x,t)}(\omega),(x,t)\in \Pi_T-\{{\rho(x,t)=0}\}.$$
 Therefore,
$$\tilde{m}^{\eps}(x,t)\rightarrow \tilde{m}(x,t), \text{a.e. in } L^{\infty}.$$

\section*{Acknowledgments}
The research is supported by China Scholarship Council No.201704910503. The author would like to heartily thank Professor Dehua Wang for helpful discussions.

\end{document}